\documentclass[10pt]{amsart}
\usepackage[]{amsmath, amsthm, amsfonts}

\newcommand {\C} {{\Bbb C}}
\newcommand {\R} {{\Bbb R}}
\newcommand {\Z} {{\Bbb Z}}

\newcommand {\HH} {{\Bbb H}}
\newcommand {\Ci} {C^\infty}

\newcommand {\cal} {\mathcal}
\newcommand {\E} {{\cal E}}

\newcommand {\I} {{\cal I}}
\newcommand {\J} {{\cal J}}
\newcommand {\K} {{\cal K}}
\newcommand {\del} {\nabla}

\newcommand {\dt} {\bullet}

\newcommand {\Om} {\Omega_{X}}
\newcommand {\Ts} {T^*}
\newcommand {\PP} {{\Bbb P}}

\newcommand {\hyp} {hyperk\"ahler }

\newtheorem{thm}[subsection]{Theorem}
\newtheorem{cor}[subsection]{Corollary}
\newtheorem{lemma}[subsection]{Lemma}
\newtheorem{prop}[subsection]{Proposition}

\newtheorem{remark}[subsection]{Remark}
\newtheorem{ex}[subsection]{Example}

\title[Geometry of cohomology support loci II]
{Geometry of cohomology support loci II:\\
 integrability of Hitchin's map} 
\author{Donu Arapura}
\address{Department of Mathematics\\
Purdue University\\
West Lafayette, IN 47907\\
U.S.A.}
\thanks{Author partially supported by the NSF}
\email{dvb@@math.purdue.edu}

\begin{document}
\maketitle

\section{ Introduction}

Green and Lazarsfeld \cite{gl2} have proven the following remarkable theorem:

\begin {thm}\label{thm:gl} Let $X$ be a smooth
complex projective variety. Then the set $S^{pq}(X)$ of line bundles in $Pic^0(X)$
satisfying $H^q(X, \Omega_X^p\otimes L) \not= 0$ is a union of a finite 
number of translates of abelian subvarieties.
\end{thm}

In this paper, we seek a generalization for higher rank bundles. The
 analogue of $Pic^0(X)$ is the moduli space $M=M_V(X,n)$ of 
semistable bundles of a given rank $n$  with trivial Chern classes, 
and within it a subset
 $\{E\,|\, H^q(X,\Om^p\otimes E)\not= 0\}$ can be defined as above.
 However, $M$ is very
far from an abelian variety in general, so it is not immediately clear what
the analogous theorem should even say. 
A clue is provided by the following theorem of Hitchin \cite{hit2}:

\begin{thm}\label{intro:thm:hit} Suppose that $X$ is a smooth 
projective curve and
 $M^s\subset  M$ the
smooth open set of stable bundles.
 Then there is morphism
from the cotangent bundle $T^*M^s$ to an affine space
 such that the general
fibers are open subsets of  abelian varieties.
\end{thm} 

Our ultimate goal then is to give  a common generalization of both theorems.
We will establish a result similar to \ref{intro:thm:hit} 
for  any component  of
 the above cohomology support locus in $M$ when
$X$ has arbitrary dimension.

In order to delve deeper into the story, it will be necessary to
explain how to compactify Hitchin's  map.  For this, we need to make the 
transition from vector bundles to Higgs bundles, which can be motivated
as follows. Suppose  that $E$ is a stable 
bundle corresponding to a smooth point $[E]\in M$. Then a cotangent
vector to $[E]$ is just a section $\theta\in H^0(X,\Omega_X^1\otimes End(E))$.
As $[E]$ is a smooth point, there is no obstruction to extending a first order
deformation of $E$ to one of  second order. The dual  
condition  is $\theta\wedge\theta = 0$. The pair $(E,\theta)$ is an 
example of a Higgs bundle.  Simpson \cite{s2} has shown that 
the set of isomorphism classes of Higgs bundles of rank $n$ with
vanishing rational Chern classes (and subject to a suitable semistability 
condition weaker than semistability of  the underlying vector bundle),
can be parameterized by a quasiprojective moduli scheme $M_{Dol}(X,n)$.
Furthermore there is a proper morphism $h$, the analogue of  Hitchin's 
map, from $M_{Dol}(X,n) $
to an affine space (which assigns to $(E,\theta)$ the characteristic
polynomial of $\theta$).  We can define
$\Sigma_{m,Dol}^k(X,n) \subseteq M_{Dol}(X,n)$ as the set of 
those pairs $(E,\theta)$ such that the appropriate $k$th cohomology
group has dimension at least $m$. In section 4, we prove the main 
result:

\begin {thm} If $X$ is a smooth projective variety, then
$\Sigma_{m,Dol}^k(X,n)$ is a Zariski closed subset
of $M_{Dol}(X,n)$. If $\tilde \Sigma$ is the normalization
of an irreducible component of $\Sigma_{m,Dol}^k(X,n)$ (with its reduced 
subscheme structure) containing a stable point, then
the connected components of the general fibers of  the pullback 
of $h$ to  $\tilde\Sigma$ are abelian varieties.
\end{thm}

The case $m=k=0$ gives an analogue of Hitchin's result. 
The above sets can be further subdivided into $(p,q)$ parts,
some components of
which form partial compactifications of the ``cotangent bundles'' of the
sets $\{E\,|\, H^q(X,\Om^p\otimes E)\not= 0\}$ considered earlier.
Similar results will be proved for these sets. 

The  analytic space associated to $M_{Dol}(X,n)$ 
has a second complex structure 
$M_B(X,n)^{an}$ which comes about via a correspondence between Higgs
bundles (of the above type) and  semisimple local systems \cite{s1,s2}. 
 The key observation is that when taken together, these  yield a  
quaternionic structure on this space, and  the set 
$\Sigma_{m,Dol}^k(X,n)$ is compatible with this structure. 
 We use this fact to show that the 
general fiber of the restriction of $h$ is lagrangian with respect to a 
suitable symplectic structure, then the theorem follows easily.
 An important precedent for the use of the 
quaternionic structure in this context
is Deligne's and Simpson's \cite{s3} approach to 
proving theorem \ref{thm:gl}. One complication, absent
 in the rank one case, is the presence of  singularities.
Recent work of Verbitsky, discussed in the  next section,
allows us to handle these issues.

 The cohomology group of a Higgs bundle has a number of different 
incarnations. To begin with, it can be defined as the hypercohomology
of an explicit complex. It is also (isomorphic to) the cohomology of the
associated local system. Both interpretations are needed in order to
verify that the cohomology support loci are quaternionic. 
The first description will also be used establish the invariance
of these loci under a natural $\C^*$-action. This will imply that any
irreducible component contains a complex variation of Hodge structure
The second
point of view will be   useful for  establishing certain homotopy
invariance properties for these sets. 
 Finally, we will reinterpret the cohomology group of a Higgs bundle 
  as an Ext group for certain sheaves on the cotangent bundle of $X$.
Then using the local to global spectral sequence, we prove a generic 
vanishing theorem in the spirit of \cite{gl1}. This leads to estimates 
on the codimension of the cohomology support loci. 

The final section of this paper contains some explicit examples. So
readers may wish to skip to it from time to time.
 For the most part,  schemes over $\C$ will be treated as sets of 
 $\C$-valued points. As usual, the superscript ``an'' indicates the 
 analytic space associated to a scheme.

\section{ Quaternionic geometry}

This section is completely expository. It is intended to give a quick
introduction to quaternionic geometry, and to some of
Verbitsky's work in particular. A nice discussion of some related ideas and
 examples
can be found in \cite{hit3}. See also \cite{hit1}, \cite{fu}, 
\cite{s4}, \cite{v1}, \cite{v2} and references contained therein.

A quaternionic (or hypercomplex) manifold is a $\Ci$-manifold
$X$ with two complex structures $\I$ and $\J$ which induce the same real
 analytic structure on $X$ and satisfies $\I\J = -\J\I$. Setting $\K = \I\J$
gives an action of the quaternions $\HH$ on the tangent bundle. Any 
quaternionic vector space is naturally a quaternionic manifold. 
A morphism of quaternionic manifolds is $\Ci$ map which is holomorphic
with respect to $\I,\J$ and $\K$ (it suffices to check holomorphicity with
respect to any two).

\begin{lemma} Let $V_1$ and $V_2$ be two finite dimensional 
quaternionic vector spaces, and let $U_i\subseteq V_i$ be open 
neighbourhoods of the origin with their induced quaternionic 
structures. Any morphism $f:U_1\to U_2$ satisfying $f(0) = 0$ is
the restriction of an $\HH$-linear map.
\end{lemma}

\begin{proof} (Deligne, see \cite{s3}). Choose a point  $x\in U_1$.
After identifying the tangent space at $x$ and $f(x)$ with $V_1$ and
$V_2$, the differential gives an $\R$-linear map $df_x:V_1\to V_2$
By assumption 
$$df_x\I = \I df_x,\quad df_x\J = \J df_x$$
so in fact $df_x$ is $\HH$-linear. Thus $df$ can be viewed as a $\Ci$ 
map from $U_1$ to $Hom_\HH(V_1,V_2)$. Differentiating again yields the
Hessian, which is a symmetric  $\R$-bilinear
form $H_x:V_1\times V_1\to V_2$.
$H_x$ is $\HH$-linear in one variable and therefore in both.
Now for the punchline:
$$\K H_x(\alpha,\beta) = \I H_x(\J\alpha,\beta)= H_x(\J\alpha,\I\beta) =
\J H_x(\alpha,\I\beta) = -\K H_x(\alpha,\beta)$$
Therefore $H_x = 0$ and the lemma follows immediately.
\end{proof}

A quaternionic submanifold of a quaternionic manifold is a 
$\Ci$ submanifold such that the inclusion is a morphism. More generally
a quaternionic subvariety $Y$ of a quaternionic manifold
$X$ is a reduced real analytic subvariety whose complexified ideal is 
locally defined by both $\I$ and $\J$ holomorphic functions.

\begin{cor} \label{cor:quat}Any quaternionic submanifold of a quaternionic vector
space $V$ is a translate of a linear subspace. A quaternionic
subvariety is a union of submanifolds.
\end{cor}

A hyperk\"ahler manifold is a $\Ci$ manifold with a Riemannian metric
$g$ and two anticommuting complex structures $\I$ and $\J$, such that 
$g$ is K\"ahler with respect to both of these structures.

\begin{prop} (\cite[6.5]{v1}) If $X$ is hyperk\"ahler then the underlying real 
analytic structures associated to $\I$ and $\J$ coincide. Therefore
$X$ is a quaternionic manifold.
\end{prop}

\begin{thm} (Verbitsky \cite{v2}). Let $(X,x)$ be a germ of a 
hyperk\"ahler manifold. The the germ of quaternionic subvariety $(Y,x)$ 
is the union of a finite number of germs of quaternionic submanifolds.
\end{thm}

Here is an outline of the proof: 
Let $R$ be the local ring of real analytic functions of $(X,x)$.
Given a complex structure $L$ on $X$, let $O_L$ be the local ring
of $L$-holomorphic functions. There is a natural inclusion
$O_L\subset R\otimes \C$ which splits: the complement is 
the ideal generated by $L$-antiholomorphic functions vanishing
at $x$. Let $\phi$ be the composition of local homomorphisms:
$$O_\I \hookrightarrow R\otimes \C \to 
O_\J\hookrightarrow R\otimes \C \to O_\I$$
A direct  calculation shows that the induced endomorphism on the
 cotangent space $m/m^2$ of $O_{\I}$
 is a homothety associated to a scalar $\lambda\in\C^*$
which is not a root of unity. After a change of variables, one can 
arrange $\I$-holomorphic coordinates so that $\phi(x_i) = \lambda x_i$.
Let $I$ be the ideal of $Y$ in $O_\I$. Then $I(R\otimes \C)$ is 
also generated by $\J$-holomorphic functions. Therefore $\phi(I) \subset I$
and this implies that $I$ is homogeneous, and thus $(Y,x)$ is 
isomorphic to the germ of its tangent cone.
The tangent cone is a quaternionic subvariety of a quaternionic 
vector space, and therefore by \ref{cor:quat} is a union of manifolds.

\begin{cor}\label{cor:norm} The normalization of a quaternionic subvariety of a 
hyperk\"ahler manifold with respect to $\I$ coincides with the 
normalization with respect to $\J$ and it  is smooth. Furthermore it 
inherits a  hyperk\"ahler structure from the ambient manifold.
\end{cor}

Given a \hyp manifold $X$, set $\omega_L(\alpha,\beta) = 
g(L\alpha,\beta)$ for $L=\I,\J,\K$. These are  just the K\"ahler forms 
associated to the complex structures.
In particular, they define (real) symplectic structures on $X$.
The form $ \omega_\J + \sqrt{-1}\omega_\K$ defines an $\I$-holomorphic
symplectic structure on $X$.

\section{Lagrangian maps}

Let $(X,\omega)$ be a real or holomorphic symplectic manifold. A closed
submanifold $Y\subset X$ is lagrangian if the tangent spaces of $Y$ are
maximal isotropic subspaces of the tangent spaces of $X$ with  respect
to the symplectic pairing induced by $\omega$. A map $h:X\to B$ 
($\Ci$ or holomorphic according to the category) to manifold $B$ will 
be called lagrangian if its differential has maximal rank, and all its
fibers are lagrangian submanifolds. Note that a lagrangian map is the
same thing as a completely integrable system. A complete discussion of
these notions would take us too far afield, see for example
\cite{gs} for further details.

If $\omega$ is only  defined on a dense open subset of $U'\subseteq X$,
a map $h:X\to B$ will be called generically lagrangian if there is a dense
nonsingular open set $U\subset B$ such that $h^{-1}(U)\subseteq U'$ and
$h^{-1}(U)\to U$ is lagrangian, the
complement of the largest such $U$ will be called the discriminant.

\begin{lemma} Let $X$ be a K\"ahler manifold with a real $\Ci$ symplectic 
structure $\omega$ (not necessarily equal to  the K\"ahler form). 
If $h:X \to B$ is a proper
lagrangian holomorphic map, then the connected components of the 
fibers are complex tori.
\end{lemma}

\begin{proof} After Stein factorization, we can assume that the fibers
are connected. Standard results in symplectic geometry \cite[page 353]{gs} show that
any fiber $F$ is diffeomorphic to a torus. It follows easily that
the Albanese map 
$$F\to H^1(F,\R)^*/H_1(F,\Z) \cong H^0(F,\Omega_F^1)^*/H_1(F,\Z)$$
is a diffeomorphism and therefore a biholomorphism. Note that the 
isomorphism $H^1(F,\R)\cong H^0(F,\Omega_F^1)$ is
only place where K\"ahler condition on $X$ is used, so the lemma
holds under  considerably weaker hypotheses.
\end{proof}

Let $(X,g,\I,\J)$ be a \hyp manifold. We will usually take $\I$ as
the preferred complex structure with  holomorphic symplectic structure
given by $\omega_\J + \sqrt{-1}\omega_K$. In particular, a lagrangian
map $h:X\to B$ will be assumed to be $\I$-holomorphic and lagrangian with 
respect to the indicated symplectic structure. Note that such a map is also 
lagrangian with respect to the real symplectic structure $\omega_J$,
consequently the above lemma applies.

\begin{thm} \label{thm:lagr}
Let $(X,g_X,\I,\J)$ and $(Y, g_Y,\I,\J)$ be a \hyp manifolds
 Suppose that $h:X\to B$ is 
 a  lagrangian map. If $f: Y\to X$ is a
  finite quaternionic morphism such that $g_Y = 
  f^*g_X$ (away from critical points). Then $h\circ f: Y \to B'$ is a  
  a generically lagrangian map, where the image $B' = h(f(Y))$ is endowed
  with the reduced analytic structure.
\end{thm}  

\begin{proof}
For any  $x\in X$, $V_x = ker\, dh_x$ is a maximal isotropic 
subspace of the tangent space $T_x$ with respect to $\omega_\J$. 
Thus $\J V_x$ is
orthogonal to $V_x$ with respect to $g_X$.
Since $$dim\, \J V_x = dim\, V_x = dim T_x/2,$$
$T_x = V_x \oplus \J V_x$.
Let $y\in Y$ be a general point, then we can identify $T_y$ with a
quaternionic subspace of $T_{f(y)}$ and $ker\, d(h\circ f)_y$ with
$V_{f(y)}\cap T_y$, and so
$T_{y} = ker\, d(h\circ f)_y\oplus \J ker\, d(h\circ f)_y$.
 Therefore $ker\, d(h\circ f)_y$ is a maximal
isotropic subspace. So the general fibers of $h\circ f$ are lagrangian.
\end{proof} 

\begin{cor} If in the above notation $h\circ f$ is proper, then
the connected components of its general fibers  are complex tori.
\end{cor}

\section{ Cohomology support loci for Higgs bundles.}

Let $X$ be a smooth complex projective variety with a fixed ample
line bundle $L$.
A Higgs bundle on $X$ consists of an algebraic vector bundle
$E$ together with a section  
$\theta\in H^0(X,\Omega_X^1\otimes End(E))$ satisfying
$\theta\wedge \theta = 0$. A  Higgs bundle $(E,\theta)$ , with 
$c_1(E)=0$ (in rational cohomology) is called stable if  
$c_1 (F).c_1(L)^{dimX -1} < 0$ for any coherent
subsheaf $F\subset E$ satisfying $rk F < rk E$ and $\theta(F)\subseteq 
\Om^1\otimes F$. A Higgs bundle with $c_1(E)= 0$ is called polystable if 
it is a direct sum of stable Higgs bundles with vanishing first Chern class.

It will be  convenient to combine the main results of Simpson 
\cite{s1, s2} into one 
big  theorem:

\begin {thm} There is an affine scheme (of finite type 
over $Spec \Z$)  $M_B(X,n)$ whose complex points
parameterize the isomorphism classes of
semisimple representation of $\pi_1(X)$ into $Gl_n(\C)$.
There is a quasiprojective scheme  $M_{Dol}(X,n)$, over $Spec \C$,
 whose complex points parameterize
polystable rank $n$ Higgs bundles with vanishing first and 
second rational Chern classes. 
There are open subsets $M_{Dol}^s(X,n)\subseteq M_{Dol}(X,n)$
and $M_{B}^{irr}(X,n)\subseteq M_{B}(X,n)$ which parameterize stable
bundles and irreducible representations respectively.
The spaces $M_{Dol}(X,n)^{an}$ and 
$M_{B}(X,n)^{an}$  are 
homeomorphic,
and $M_{Dol}^s(X,n)$ and $M_{B}^{irr}(X,n)$ correspond under this 
homeomorphism.
\end{thm}

\begin{remark} 
0) $X$ and $n$ will be omitted from the notation when it is safe
to do so.

1) The word ``parameterize'' is a bit vague. The
correct statement is that these are coarse moduli spaces for the
appropriate moduli functors.

2) A semistable Higgs bundle with $c_1=0$
is an iterated extension of stable Higgs bundles with $c_1=0$. 
Two semistable bundles are equivalent if the their stable
factors coincide (up to isomorphism). Every equivalence
class of semistable bundles has a unique polystable representative.
Thus $M_{Dol}$ parameterizes equivalence classes of semistable
bundles. Similarly $M_B$ parameterizes equivalence classes of arbitrary
representations, where two representations are equivalent if they
have isomorphic semisimplifications.

3) These moduli spaces may be nonreduced. However we will usually suppress the 
scheme structure and just treat them as sets of $\C$-valued points.
 For every (poly, semi)stable Higgs bundle
or semisimple representation $V$, let $[V]$ denote the corresponding
point in the moduli space.

\end{remark}

A family of Higgs bundles on $X$  parameterized by $T$ is a 
vector bundle ${\cal E}$ on $X\times T$,
with a  section $\Theta$ of $p_X^*\Om^1 \otimes End({\cal E})$
satisfying $\Theta\wedge\Theta = 0$.
$M_{Dol}$ is only a coarse moduli space, so there  may not be
a universal family of Higgs bundles. However Simpson's construction gives a bit
more. Namely $M_{Dol}(X,n)$ is a quotient, in the sense of geometric invariant 
theory, of a
locally closed subscheme $Q_n$ of an appropriate Quot or Hilbert scheme.
$X\times Q_n$ will in fact carry a  family of Higgs bundles
$({\cal E}, \Theta)$ 
such that  its restriction  $({\cal E},\Theta)|_q$
to a slice $X\times \{q\}$ is semistable and corresponds
to the image of $q$ in $M_{Dol}(X,n)$. Over $M_{Dol}^s(X,n)$ it 
possible to find 
cross sections to $Q_n\to M_{Dol}(X,n)$ locally in the etale topology.
The following is an almost immediate consequence:

\begin{prop}\label{prop:closed} 
Let $S$ be a set of semistable Higgs bundles on $X$.
Suppose that for every $T$ and semistable family of
Higgs bundles $({\cal E},\Theta)$
parameterized by $T$, the set
 $S_T = \{t\,|\, ({\cal E},\Theta)|_t \in S\}$ is 
Zariski closed. Then the set $S_M$ of all $m\in M_{Dol}$ which possess a
 semistable representative in $S$ is Zariski closed.
\end{prop}

\begin{proof} $S_M$ is constructible because it
is the image of $S_Q$ under the canonical map. Given a limit point 
$s$ of $S_M$, there is an irreducible curve $C$ such that $s\in C$ and 
$C-\{s\}\subset S_{M}$.
As $Q\to M_{Dol}$ is surjective, there is an irreducible curve $C'\subset Q$ and
finite map $C' \to C$. One obtains a family of Higgs bundles
on $C'$ by restriction. $S_{C'}= C'$ since it is closed and contains 
the preimage of $C-\{s\}$. Therefore $s\in S_{M}$ and so $S_M$ is closed.
\end{proof}

\begin{prop}\label{prop:restrict}
Let $C\subset X$ be a curve 
obtained as a complete intersection
of divisors associated to an $N$th power of $L$, with $N >> 0$.
Then the restriction maps
$M_B(X,n) \to M_B(C,n)$ and $M_{Dol}(X,n) \to M_{Dol}(C,n)$
are injective (on $\C$-valued points)
and compatible with the homeomorphisms of the associated analytic spaces.
\end{prop}

\begin{proof} This well known to the experts, so we
will merely indicate the main ideas.
The first part follows from the Lefschetz hyperplane
theorem \cite{mi}. For the second part, first note that the restriction
of a semistable Higgs  bundle is semistable thanks to Simpson's
generalization of the Mehta-Ramanathan theorem \cite{s1}, 
so the map is well defined. 
As for injectivity, choosing $N >>0$ guarantees that 
$$H^1(E_1^*\otimes E_2\otimes I_C) = 0$$
for any pair of Higgs bundles on $X$. Thus any isomorphism of their 
restrictions can be lifted to map of the Higgs bundles 
which can be seen to be
an isomorphism (using, for example, the fact that polystable bundles 
are direct sums of simple bundles). The last part is just a 
restatement of the naturality of the correspondence.
\end{proof}

\begin{remark} The images of $M_{Dol}^s(X,n)$ and 
$M_B^{irr}(X,n)$ under restriction lie in the corresponding subsets
for $C$.
\end{remark}

Given a Higgs bundle $(E,\theta)$, define $H^i(E,\theta)$ to be
the $i$th hypercohomology of the 
 complex:
$$
E\stackrel{\theta\wedge}{\longrightarrow}
\Omega_X^1\otimes E \stackrel{\theta\wedge}{\longrightarrow} 
\Omega_X^2\otimes E \ldots 
$$

We define the cohomology support loci as:
$$\Sigma^k_{m,Dol}(X,n) = 
\{[(E,\theta)] \in M_{Dol}(X,n)\,|\, 
(E,\theta)\>{ polystable\> and\>}dim\, H^k(E,\theta) \ge m\} $$
We can make an analogous definition for $M_B$. Given a representation
of $\rho:\pi(X)\to Gl_n(\C)$, there exists (up to isomorphism) a unique
rank $n$ locally constant sheaf, or local system, with $\rho$ as
its monodromy representation. So we can, and will, view $M_B(X,n)$ as a moduli
space of local systems on $X$. Let
$$\Sigma^k_{m,B}(X,n) =
\{ V \in M_B(X,n)\, | \, 
V\> { semisimple\> and}\> dim\, H^k(X,V) \ge m\}$$

\begin{prop} $\Sigma^k_{m,Dol}$ is Zariski closed in
$M_{Dol}$. $\Sigma^k_{m,B}$ is is Zariski closed in
$M_B$. $\Sigma^k_{m,Dol}$ and $\Sigma^k_{m,B}$ coincide under
the correspondence between $M_{Dol}^{an}$ and $M_B^{an}$.
\end{prop}

\begin{proof} 
Note that if $(E_i,\theta_i),\, i=1,2$ are equivalent semistable
bundles with $(E_1,\theta_1)$ polystable, then
$$dim H^k(E_1,\theta_1) \ge dim H^k(E_2,\theta_2)$$
by subadditivity of cohomology. 
Therefore the first  statement  follows from proposition 
\ref{prop:closed} and 
 the  semicontinuity theorem for cohomology \cite[7.7.5, 7.7.12]{ega3}
 The second statement is proved in \cite{a2}. The last part follows
 from \cite[2.2]{s1}.
\end{proof}

\begin{thm}(Hitchin\cite{hit1}) Let $C$ be a smooth projective curve.
Then the stable locus  $M_{Dol}^{s,an}(C,n) $ is smooth and 
carries a natural hyperk\"ahler structure $(g,\I,\J)$. Where $\I$ is the usual
complex structure and $\J$ is the structure induced from the 
identification $M_{Dol}^{s,an} \cong M_B^{irr,an}$.
\end{thm}

\begin{remark} Hitchin  stated this only  when $n=2$, although his 
proof presumably works for higher rank bundles. 
In any case, Fujiki \cite{fu} has established a much more general 
result.
\end{remark}

Given a rank $n$ Higgs bundle $(E,\theta)$ on $X$, let
$\theta^{i}\in H^{0}(X, S^{i}\Om^{1})$ be the $ith$ (symmetric) power
Let
$$h_X(E,\theta)\in S_n(X) = \bigoplus_{i=1}^n\, H^0(X,S^i\Omega_X^1)$$
be the map $ (E,\theta) \mapsto trace(\theta^{i})$.
This differs from Simpson's and Hitchin's 
definition in that they use the 
characteristic polynomial. However there is 
an algebraic automorphism $\sigma$ of the target space such that 
$\sigma\circ h$ agrees with their map, so they are essentially the same.
Therefore:

\begin{thm} (Simpson\cite{s2}, Hitchin\cite{hit2}) $h_X$ gives a
proper morphism from $M_{Dol} $
to $S_n(X)$. When $dim X = 1$, $h_X$ is surjective and
generically lagrangian with respect to the
 hyperk\"ahler structure on  $M_{Dol}^{s,an}$.
\end{thm}

Suppose that a curve $C\subseteq X$ has been chosen as in
\ref{prop:restrict},
then we obtain a commutative diagram:

$$\begin{array}{ccc}
M_{Dol}(X,n)  & \hookrightarrow   & M_{Dol}(C,n)\\
h_X\downarrow &                   & h_C\downarrow \\
S_n(X)        & \stackrel{i}{\to} & S_n(C) \\
\end{array} $$
While $i$ need not be injective, we do have:

\begin{lemma} The map $h_X(M_{Dol}(X,n))\to h_C(M_{Dol}(C,n))$ is
finite.
\end{lemma}

\begin{proof} This follows from the properness of all the other maps
in the diagram. 
\end{proof}

Consequently, fibers of  $h_X$ are components of fibers of $i\circ h_X$.
Putting the above results together yields the main theorem:

\begin {thm} Let $X$ be a smooth projective variety.
 If $\tilde \Sigma $ is a connected
 component of the normalization
of $\Sigma^k_{m,Dol}$ (with its reduced subscheme structure) which 
meets $M_{Dol}^s$,
and $\tilde h:\tilde \Sigma \to h(\tilde \Sigma)$ the
natural map.
Then $\tilde h$ is generically lagrangian.
 In particular the connected 
components of its general fibers 
 are abelian varieties of dimension half of that of $dim\, 
 \tilde\Sigma$
\end{thm}

\begin{proof} Let $C\subseteq X$ be a complete intersection curve of high
degree. Then $\Sigma_{m,Dol}^k(X,n)$ is a quaternionic subvariety of
$M_{Dol}(C,n)$ by the previous results. Therefore the  theorem follows
from \ref{thm:lagr}.
\end{proof}

\begin{cor} \label{cor:fibdim} Any  irreducible component $F$ of every
fiber  of $\tilde h$ satisfies $dim F = dim\, \tilde \Sigma/2$.
\end{cor}

\begin {proof} By upper semicontinuity of dimensions, it is enough to
prove $dim F \le  dim\, \tilde \Sigma/2$. By \ref{cor:norm}, $\tilde 
\Sigma$ is compatible with the complex structure coming from
$\Sigma_{m,B}^k(X,n)$. With this structure $\tilde \Sigma^{an}$ is Stein,
therefore $H_i(\tilde \Sigma^{an}, \Z) = 0$ for $i > dim \tilde S$ \cite{n}.
On the other hand, with the original complex structure 
$\tilde \Sigma$ is quasiprojective, and $F$ is a proper subvariety.
 Therefore the fundamental class of $F$ defines a nonzero element of
 $H_{2dim(F)}(\tilde\Sigma^{an},\Z)$, and this forces the inequality.
\end{proof}

\begin{remark}  The theorem can also be used to recover a result of Biswas 
\cite[8.2]{b} about Poisson commutivity of higher dimensional Hitchin
maps.
\end{remark}

In order to get a better analogue of the original theorem of Green and
Lazarsfeld, we need to break $\Sigma_m^k$ up into $(p,q)$ parts.
For a Higgs bundle $(E,\theta)$, let
$H^{pq}(E,\theta)$ be the $E^{pq}_\infty$ term
of the spectral sequence:
$$E_1^{pq} = H^q(X,\Omega_X^p\otimes E) \Rightarrow H^{pq}(E,\theta)$$
Then set $S^{pq}_m(X,n)$ equal to the closure of
$$\{(E,\theta) \in M_{Dol}(X,n) | dim\, H^{pq}(E,\theta) \ge m\}$$
Clearly
$$\Sigma_{m,Dol}^k(X,n) 
 = \bigcup_{ m_0+m_1+\ldots = m}\bigcap_{p}\, S^{p,k-p}_{m_p}(X,n)$$

\begin{cor} Let  $\tilde S$ be  a connected component of the 
normalization of $S_m^{pq}$ which meets $M_{Dol}^s$.
The connected components of
general fibers of the restriction of $h$ to  $\tilde S$ 
are abelian varieties.  All fibers have dimension equal to one half 
of that of $\tilde S$.
\end{cor}

This can be deduced from the theorem and the previous corollary
 using the next lemma:

\begin{lemma} Let $X$ be a noetherian topological space. Suppose that
there are nested closed sets 
$$X = X^i_0 \supseteq X^i_1 \supseteq \ldots$$
$i = 0,\ldots n$. Then any irreducible component of $X^i_m$ is an 
irreducible component of some set of the form
$$Y_m  = \bigcup_{m_0+m_1+\ldots = m}\bigcap_{ i}\,X_{m_i}^i$$
\end{lemma}

\begin{proof} Suppose that $S$ is an irreducible component of $X_{m_0}^0$.
We can assume that $S$ is not contained in $X_{m_0+1}^0$. Similarly, 
let $m_1,\ldots$ be the largest integers for which $S$ is contained 
in $X_{m_1}^1,\ldots$. Then is easily seen to be an irreducible component of 
$Y_m$, where $m = m_0+m_1+\ldots$.
\end{proof}

If $(E,\theta)$ is a Higgs bundle, its dual is $(E^*,-\theta)$ where
$\theta$ is viewed as section of $\Om^1\otimes End(E^*)\cong
\Om^1\otimes End(E)$. This is compatible with the duality of local
systems. Simpson has already observed a duality theorem holds on
cohomology. This can be refined slightly:

\begin{prop} If $d = dim\, X$ then
$$H^{pq}(E,\theta) \cong  H^{d-p,d-q}(E^*,-\theta)^*$$
\end{prop}

\begin{proof} Set $\omega_X = \Omega_X^d$.
A special case of the Grothendieck-Serre duality theorem is that
$$\HH^{d-i}({\cal H}om(V^\dt,\omega_X)) \cong \HH^i(V^\dt)^*$$
for any finite complex of locally free sheaves $V^\dt$.
The natural pairing of $\Om^p\otimes \Om^{d-p}\to \omega_X$  induces an
 isomorphism of complexes
$${\cal H}om((\Om^\dt\otimes E,\theta), \omega_X) \cong 
(\Om^\dt\otimes E^*,-\theta)[d]$$
upto sign. This  respects the the filtrations
$${\cal H}om((\Om^{\le p}\otimes E,\theta), \omega_X) \cong 
(\Om^{\ge d-p}\otimes E^*,-\theta)[d]$$
and induces isomorphisms on the associated graded parts. Therefore
there is an isomorphism between the associated spectral sequences
converging to the hypercohomology groups on the left and right hand
sides. The $E_{\infty}$ terms can be identified with the groups of the 
proposition. 

\end{proof}

\begin{cor}\label{cor:duality} $S^{pq}_m(X,n) = S^{d-p,d-q}_m(X,n)$
\end{cor}

\section{ Cohomology support loci for vector bundles.}

We will use the same notation as in the previous section. Let $M_V(X,n)$ 
be the closed subscheme of $M_{Dol}(X,n)$ parameterizing polystable
Higgs bundles of the form $(E,0)$. 
Let $M^s_V$ be the open subset of stable bundles.
Note that $(E,0)$ is (poly, semi)stable
if and only if $E$ is (poly, semi)stable in the usual 
sense (with respect to ``slope'').
 So $M_V$ is just the coarse moduli space of semistable
vector bundles of rank $n$ with trivial Chern classes.
Note that $(E,\theta)$ is (poly, semi)stable if $E$ is.
Let  $T^*M_V$ (respectively $T^*M_V^s$) be the set of all
 Higgs bundles $[(E,\theta)]$ such that
$E$ is polystable (respectively stable). 
 $M_V$ may have singularities, so the notation $T^*M_V$ is 
merely suggestive; some justification for it will be given below.
There is a morphism $T^*M_V \to M_V$ given
by projection.  

A vector bundle is determined by a $1$-cocycle $g_{ij}\in Z^1({\cal U}, Gl_{n}(O_{X}))$,
and a  first order deformation to it is described by a cocycle of the form
$$g_{ij}+\epsilon\gamma_{ij}\in Z^1({\cal U},Gl_{n}(O_{X}[\epsilon]/(\epsilon^{2})))$$
$\gamma_{ij}$ defines a cocycle with values in $End(E)$.
So the Zariski tangent space to $[E] \in M_{V}^s$ can be identified 
with $H^{1}(X,End(E))$.
The obstruction to lifting a first order deformation $v\in H^{1}(X,End(E))$ to 
one of second order is given by 
$[v, v]\in H^{2}(X,End(E))$.
 There are in fact no higher obstructions, 
so the tangent cone to $[E]\in M_{V}^s$ is just 
$\{v\in H^{1}(X,End(E))\, |\, [v, v] = 0\}$ \cite[9.4]{gm}.
By a theorem of  Donaldson \cite{d} and Uhlenbeck-Yau \cite{uy}, $E$ 
carries unitary flat connection $\del$. Consequently, $End(E)$ also 
carries a unitary flat connection. Hodge theory with unitary flat 
coefficents and the self duality of $End(E)$ shows that there is a 
conjugate linear isomorphism
$$H^{i}(X,End(E)) \cong H^{0}(X,\Om^{i}\otimes End(E))$$
preserving the graded Lie brackets. The Higgs fields on $E$ are
precisely the conjugates of the vectors in the tangent cone.
This implies that they have the same dimension as real algebraic 
varieties, and therefore as complex algebraic varieties.
It is sometimes better to view
$T^{*}_{[E]}=H^{0}(X,\Om^{1}\otimes End(E))$ as dual to the tangent space,
via the hard Lefschetz pairing 
$$<\alpha,\beta> = \int_{X}trace(\alpha\cup \beta)\cup L^{dimX -1}$$
In an analogous fashion , the tangent space to any stable point
$[(E,\theta)]\in M_{Dol}$ is $H^{1}(End(E,\theta))$ where 
$End(E,\theta)$ is the Higgs bundle $(End(E),ad(\theta))$ \cite[10.5]{s2}.
The tangent cone is defined by the quadratic form associated to the
cup product
$$H^{1}(End(E,\theta)) \times H^{1}(End(E,\theta)) \to 
H^{2}(End(E,\theta)).$$

\begin{prop}\label{prop:cotang} $T^*M^s_V$  is an open subset of $M_{Dol}$.
If $M\subseteq M_{V}$ is an irreducible component, then $dim M = 
dim (\pi^{-1}M)/2 $ and $M$ is an irreducible component of
$h^{-1}(0)\cap \pi^{-1}M$.
\end{prop}

\begin{proof} 
As pointed out in the remarks preceding \ref{prop:closed},
etale local cross sections to $Q^s\to M_{Dol}^s$ exist. Therefore
 there is an 
etale neighbourhood $T \to M_{Dol}$ of any stable point $(E_{0},\theta_{0})$
 and family of Higgs bundles 
$(E_{t},\theta_{t})$ parameterized by $T$, such that 
$[(E_{t},\theta_{t})]$ is precisely the image of $t$. Stability is an 
open condition \cite[3.7]{s2},
thus if  $E_{0}$ were stable, then this would hold in open 
neighbourhood of $0\in T$.

Note that $M_{V} \subseteq h^{-1}(0)$, and the dimension of any component 
of $h^{-1}(0)$ is half the dimension of an irreducible component of 
$M_{Dol}$ containing it by \ref{cor:fibdim}. Thus it suffices to prove that 
$dim M \ge 2 dim(\pi^{-1}M)$.
Choose a general point $[E]\in M$. Then consider the terms of low 
degree for the spectral sequence converging to the cohomology of 
$End(E,0)$:
$$0\to H^{0}(X,\Om^{1}\otimes End(E)) \to H^{1}(End(E,0)) \to 
H^{1}(X, End(E)) \to 0$$
It follows that the tangent cone $C_{1}$ of $M_{Dol}$ at $[(E,0)]$ maps to the 
tangent cone $C_{2}$ of $M_{V}$ at $[E]$, and the fiber over $0$ is the 
cone $\Ts_{[E]}$. This implies that the dimension of $C_{1}$ is less than 
or equal to $dim C_{2} + dim \Ts_{[E]} = 2dim C_{2}$.
\end{proof}

The intersections of $S^{pq}_m$ with $M_V$ and
$T^*M$ have a rather concrete description which is very close
to the spirit of \cite{gl1}. Let us write $T^*S^{pq}_{m}$ for
$S^{pq}_m \cap T^*M_V$.

\begin{prop} If $(E,\theta)$ is polystable, then $[(E,\theta)] \in T^*S^{pq}_{m}$
if and only if the $p$th cohomology of the complex
$$\ldots H^q(X,\Om^p\otimes E)\stackrel{\theta\wedge}{\longrightarrow}
H^q(X,\Om^{p+1}\otimes E)\stackrel{\theta\wedge}{\longrightarrow} \ldots$$
has dimension greater than or equal to  $m$. In particular
$(E,0)\in S^{pq}_1$ if and only if $H^q(X,\Om^p\otimes E) \not=0$.
\end{prop}

The proposition is a consequence of the next two lemmas.

\begin{lemma} If $E$ is polystable, then $H^{pq}(E,\theta)$ is just 
the $p$th cohomology of the complex
$$\ldots H^q(X,\Om^p\otimes E)\stackrel{\theta\wedge}{\longrightarrow}
H^q(X,\Om^{p+1}\otimes E)\stackrel{\theta\wedge}{\longrightarrow} \ldots$$
\end{lemma}
\begin{proof} 
This is equivalent to the assertion that the spectral sequence
$$E^{pq}_1 = H^q(X,\Om ^p\otimes E) \Rightarrow H^{p+q}(E,\theta)$$
degenerates at $E_2$.
By the theorem of Donaldson  and Uhlenbeck-Yau,
$E$ carries a unitary flat connection $\nabla$. The lemma now follows by
applying \cite[III 3.6]{a2} to the complex constructed in the proof
of [loc.cit, IV 2.1]. (In the notation of that paper, the spectral
sequence associated to ${\cal V}(0)$ coincides with the one above.)

An alternative argument can be given by modifying the proof of 
\cite[3.7]{gl1} by replacing $\partial$ and $\bar\partial$ by the
$(1,0)$ and $(0,1)$ parts of $\nabla$.
\end{proof}

\begin{lemma} If  $(E_t,\theta_t)$ is a family of Higgs bundles
parameterized by a smooth curve $T$ such that $E_0$ is
polystable for some $0\in T$ and $dim \,H^{pq}(E_t,\theta_t) \ge m$
for $t\not= 0$ then $dim\, H^{pq}(E_0,\theta_0) \ge m$.
\end{lemma}

\begin{proof}
We can assume that $T = Spec R$. 
Consider the complex of $R$-modules
$$\ldots H^q(X,\Om^p\otimes E)\stackrel{\theta\wedge}{\longrightarrow}
H^q(X,\Om^{p+1}\otimes E)\stackrel{\theta\wedge}{\longrightarrow} \ldots$$
Our  assumptions imply that if one tensors this by the residue field
of $t\not= 0$, the $p$th cohomology has dimension greater than or equal
to $m$ (This is true regardless of whether the spectral sequence 
degenerates for $t\not= 0$, because at any rate $dim\, E_2\ge dim 
E_\infty$). Therefore this property persists for $t=0$, and the lemma
follows from the previous one.
\end{proof}

The next result gives a useful dimension estimate on the cohomology 
support loci. It can also be deduced using Green's and Lazarsfeld's  
deformation theory \cite{gl1}.

\begin{prop}\label{prop:dimest} Let  $S_V$  be an irreducible
component of $S^{pq}_{m}\cap M_V$, and let $S$ be the irreducible component of $S^{pq}_m$ 
containing $S_V$.
Then for a general point $[E] \in  S_V$ 
$$dim S_V = dim (T^*_{[E]}\cap S) =
 \frac {dim S}{2}$$
\end{prop}

\begin{proof} By \ref{cor:fibdim} and \ref{prop:cotang},
$dim S_V = \frac {dim S}{2}$. The remaining equality follows
from $dim S = dim S_V + dim (T^*_{[E]}\cap S) $.
\end{proof} 

\section{ $\C^*$-invariance}

In the last section, we made an explicit study of the intersection of
the cohomology support loci with $T^*M_V^s$. There are some features
of this geometry which extend  to the whole space. First, recall:

\begin{thm} (Simpson \cite{s2}) There is an algebraic $\C^*$-action
on $M_{Dol}$ given by $t:(E\theta) \mapsto (E,t\theta)$. For any
point $e\in M_{Dol}$ the limit of $ te$ as $t\to 0$ exists.
\end{thm}

Of particular interest are the fixed points. A Higgs bundle
$(E,\theta)$ is called a complex variation of
Hodge structure if $E$ admits a grading $\oplus E^p$ such that 
``Griffiths transversality'' $\theta(E^p) \subset E^{p-1}$
holds. Given a complex variation of Hodge structure,
let $T$ be the automorphism which acts by $t^{-p}$ on $E^p$ where
$t\in \C^*$. Then
$T$ induces an isomorphism  $(E,\theta)\cong (E,t\theta)$. Therefore
$[(E,\theta)]$ is a fixed point. Conversely, Simpson has shown that all
fixed points arise this way, and in fact if the underlying Higgs bundle
$(E,\theta)$ is stable then the grading is uniquely determined (up to
to a shift of indices).
Note that if $\theta = 0$, then we can take $E^0 = E$

The previous theorem implies that the $\C^*$-action extends to a
morphism of reduced schemes
${\Bbb A}^1\times M_{Dol,red} \to M_{Dol,red}$. The image of
$\{0\}\times M_{Dol}$ is the fixed point set $F$. Thus we
get a morphism $\pi:M_{Dol,red} \to F_{red}$ by composing
$$M_{Dol,red} \cong \{0\}\times M_{Dol,red} \to F_{red}.$$
$\pi$ extends the map $T^*M_V\to M_V$ constructed earlier,
and exhibits $M_{Dol}$ as a family of cones.
The cohomology support loci are compatible with this conical
structure:

\begin{lemma} $\Sigma^k_{m,Dol}$ is invariant under the $\C^*$-action.
\end{lemma}
\begin{proof} 
There is an isomorphism $H^k(E,\theta)\cong H^k(E,t\theta)$
given on the level of complexes by 
multiplication by $t^{p}$ on  $\Om^p\otimes E$.
\end{proof}

As these sets are closed, we obtain:

\begin{cor} $\pi(\Sigma^k_{m,Dol}) \subset \Sigma^k_{m,Dol}$, so 
$\Sigma_{m,Dol}^k$ is 
a family of subcones.
\end{cor}

\begin{cor} Any irreducible component of $\Sigma^k_{m,Dol}$ contains a
complex variation of Hodge structure.
\end{cor}

Conjectures of Simpson \cite{s1} and Pantev \cite{p} suggest that much more
should be true, for example every component should contain an integral
variation of Hodge structure.
In the rank one case, the $\C^*$ invariance of the above sets is a
powerful constraint. In fact, it leads to a proof
of theorem \ref{thm:gl}, \cite{a1,s3}.

\section{ Generic vanishing}

 In this section, we relax the condition on the Higgs bundles. 
We no longer
insist that the Chern classes vanish, or even that they are locally 
free.
Let $X$ be a smooth projective $d$-dimensional variety
A Higgs sheaf on $X$ is a torsion free coherent $O_{X}$-module $E$ 
together with morphism $\theta:E \to \Om^{1}\otimes E$ satisfying
$\theta\wedge\theta = 0$. There is a useful alternative viewpoint 
which we now recall.
Let  $\pi: 
T^{*}X \to X$ be the cotangent 
bundle. A cotangent vector $\eta\in T^{*}X$ can also be viewed as
an element of the fiber of $\pi^*\Om^{1}$.
This defines  a canonical section $\Theta$ of $\pi^*\Om^{1}$. If 
$\cal E$ is a torsion free $O_{T^{*}X}$-module, such that $supp(\E)\to 
X$ is finite, then $E = \pi_{*}{\cal E}$ is a torsion free coherent 
sheaf, and the map
$$ [E\stackrel{\theta}{\rightarrow} \Om^1 \otimes E ]
\cong \pi_{*}(\E \otimes [O_{T^{*}X}\to \pi^*\Om^{1}])$$
defines a Higgs structure on $E$. In fact, every Higgs sheaf arises 
this way from a unique $\E$ \cite{s1,s2}. The 
support of $\E$ is 
precisely the set of eigenforms for $\theta$. In other words,
$supp(\E)\cap T^{*}_{x}$ is the set
 of cotangent vectors $\eta$ satisfying $\theta_{x}(v)= \eta v$ 
for some nonzero $v\in E_{x}$. Define the degeneracy locus of 
$(E,\theta)$ as $supp( \E)\cap X$ where
$X\subset T^{*}X$ is identified with the zero section. More concretely,
if $E$ is locally free then
$x$ lies in the degeneracy locus if and only if $\theta_{x}$ has a zero eigenform or,
equivalently, is not injective. Let $degen(\theta)$ be the dimension 
of the degeneracy locus of $(E,\theta)$ if it is nonempty, or $-1$
otherwise. It will simplify matters to  define $dim\emptyset = -1$.

Given a Higgs sheaf $(E,\theta)$, the cohomology $H^{i}(E,\theta)$ can 
be defined as the hypercohomology of the complex 
$$
E\stackrel{\theta\wedge}{\longrightarrow}
\Omega_X^1\otimes E \stackrel{\theta\wedge}{\longrightarrow} 
\Omega_X^2\otimes E \ldots 
$$
as  before.

\begin{prop}\label{prop:ext} If $X\subset T^{*}X$ is identified with the zero 
section, then for any Higgs sheaf $(E,\theta)$ and corresponding 
sheaf $\E$ on $\Ts X$,
$$H^{i}(E,\theta) \cong Ext^{i}(O_{X},\E)$$
\end{prop}

\begin{proof}
The zero locus of $\Theta$ is precisely $X$. Therefore $O_{X}$ is 
quasiisomorphic to the Koszul complex
$$K_{\bullet} = \ldots \wedge^{2}\pi^{*}T_{X}\to \pi^{*}T_{X} \to O_{T^{*}X} \to 0$$
Therefore 
$$Ext^{i}(O_{X},\E) \cong H^{i}({\cal H}om(K_{\dt}, \E))
\cong H^{i}(\pi_{*}(K_{\bullet}^{*}\otimes \E))$$
By the projection formula,
$$\pi_{*}(K_{\bullet}^{*}\otimes \E) = \Omega_{X}^{\bullet}\otimes E$$
and  the differentials are easily seen to be given by $\theta\wedge$.
\end{proof}

\begin{thm} If $(E,\theta)$ is a Higgs sheaf and $i > degen(\theta) + 
d$, then $H^i(E,\theta) = 0$. If $E$ is locally free then
$H^i(E,\theta)=0$
for $i < d -degen(\theta)$ also.
\end{thm}

\begin{proof} The support of ${\cal E}xt^{\bullet}(O_X,\E)$ lies 
in  the degeneracy locus $supp(\E)\cap X$.
 Therefore the first part of the theorem follows from \ref{prop:ext}
and the spectral sequence
$$ H^p({\cal E}xt^q(O_X,\E)) \Rightarrow Ext^{p+q}(O_X,\E).$$

The second statement follows by duality \cite[2.5]{s1}.
\end{proof}

As an immediate corollary: $\Sigma^{i}_{1} \not= M_{Dol}$ if there 
exists a Higgs bundle $(E,\theta)$ with $  |i -d| > degen(\theta) $.
When $E$ is stable, this can be interpreted as a generic vanishing
theorem:

\begin{cor} If $[E]\in M^{s}_{V}$, and $\theta$ is a Higgs field on 
$E$. Then in any neigbourhood of $[E]$ in $M^{s}_{V}$ there is an $[E']$ 
such that 
$$H^{q}(X,\Om^{p}\otimes E') = 0$$
for $|p + q -d|> degen(\theta)$.
\end{cor}

One gets a very concrete statement by restricting to diagonal 
Higgs fields.

\begin{cor} If $\phi \in H^{0}(X,\Om^{1})$, then for any irreducible
component  $M$ of $M_{V}(X,n)$, there exists a semistable vector
bundle $E$, with $[E]\in M$ and
$$H^{q}(X,\Om^{p}\otimes E) = 0$$
for $|p + q-d| >  dim \{x\,|\, \phi_{x} = 0\} $.
\end{cor}

\begin{proof} Apply the theorem to  $(E,\theta)$ where  $[E]\in M$ is a 
general point and  $\theta = \phi I_{E}$.
\end{proof}

The above arguments can be refined to yield a codimension estimate. 
Given a Higgs field $\theta$, clearly $degen(\lambda\theta) = 
degen(\theta)$ for any  $\lambda \not= 0$. Thus $degen$ gives a map,
from the projectivization $\PP\Ts_{[E]}$ to the set of natural numbers,
which is easily seen to be upper semicontinuous.

\begin{cor} If $S \subseteq \Sigma^{i}_{1}\cap M_{V}^s$ is an irreducible 
component, and $[E]\in S$ a general point. Assume that $M_{V}$ is 
smooth at $[E]$, then the codimension of $S$, 
in the irreducible component of $M_{V}$ containing it, is greater than
  $dim\, \{[\theta] \in \PP\Ts_{[E]}\, | \, degen(\theta) \ge |i-d| \}.$
\end{cor}

\begin{proof} The smoothness assumption implies that $\Ts_{[E]}$ is a vector space.
By \ref{prop:dimest}, 
 the codimension of the projectivized cone $V = \PP(\Ts_{[E]}\cap 
\Sigma^{i}_{1})\subset \PP \Ts_{[E]}$ coincides with the above 
codimension.
By the theorem,  $V$ cannot meet $\{[\theta] \in \PP\Ts_{[E]}\, | \, 
degen(\theta) \ge |i-d| \}$. Thus the corollary is a consequence of 
Bezout's theorem.
\end{proof}

\section{Homotopy invariance}

The cohomology support loci are clearly isomorphism invariants.
But much more is true, namely the $\Sigma_{m,B}^k$
depend only on the homotopy type of the space. We will give some
variants of this which will be quite useful for the construction of
examples.  

\begin{prop} Suppose that $f:X \to Y$ is a morphism of smooth projective 
varieties, with connected fibers, such that the induced map on homotopy groups
$\pi_i(X^{an})\to \pi_i(Y^{an})$ is an isomorphism  for $i \le k$
and a surjection for $i = k+1$ with $k >0$. Then  for all $m$,
 $f^*(\Sigma_{m,Dol}^i(Y))$ is contained in $ \Sigma_{m,Dol}^i(X)$
when $i = k+1$ and equality holds when $i \le k$.

\end{prop}

\begin{proof} It suffices to prove the corresponding statement for
$\Sigma_B$. First note that by standard arguments in topology
\cite[pages 99-100]{sp} $f^{an}$ is homotopy 
equivalent to a fibration of topological spaces $f':X'\to Y'$.
The homotopy long exact sequence  and the  hypothesis
implies that 
the first $k$  homotopy groups of the fiber $F$
vanish. Therefore by Hurwitz theorem $H^j(F,V) = 0$ for $j \le k$
and $V$ and arbitrary coefficient group. Thus the  Leray-Serre spectral 
sequence yields an injection
$$H^i(Y^{an},V) \to H^i(X^{an},f^*V)$$
when $i = k+1$ and an isomorphism when $i \le k$.
\end{proof}

\begin{cor}\label{cor:lef} If $X\subset Y$ is general hyperplane
section where $dim Y \ge 3$, then under restriction:
$\Sigma_{m,Dol}^k(Y) \subseteq \Sigma_{m,Dol}^k(X)$ for $k =
dim X$ and equality holds when $k < dim X$.
\end{cor}

\begin{proof} The hypothesis of the theorem is fulfilled by the Lefschetz hyperplane
theorem \cite{mi}. 
\end{proof}

\begin{prop} If $f:X\to Y$ is a surjective map of smooth projective 
varieties, then $f^*\Sigma_{m,Dol}^k(Y) \subseteq \Sigma_{m,Dol}^k(X)$ for
$k \le dim Y$.
\end{prop}

\begin{proof} Once again, we will work with $\Sigma_B$. We will also 
omit the superscript ``an''.
 We will show that the map $f^*:H^k(Y,V)\to H^k(X,V)$ is injective 
for any local system $V$.
To begin with, assume that $f$ is generically finite.
Then one gets a transfer or Gysin homomorphism $f_!:H^k(X,V)\to 
H^k(Y,V)$ as the Poincar\'e dual of $H^{2dimY - i}(Y,V^*)\to 
H^{2dimY-i}(X,V^*)$. Arguing as in the case of constant coefficients,
$\frac{1}{deg\,f}f_!$ splits $f^*$.

In general, choose a complete intersection $Z\subset X$ of hyperplane 
sections, such that $f|_Z$ is finite. The map
$H^k(Y,V)\to H^k(Z,f^*V)$ factors through $f^*$. Consequently $f^*$ is
again injective.
\end{proof}

We  extend the notation  $\Sigma_B^k(X)$ so as to
allow $X$ to be any topological space. Of course, there won't be an
analogue $\Sigma_{Dol}$ in general. Our next task is to prove
a K\"unneth decomposition. Given two spaces $X_1, X_2$, there is a
morphism 
$$\tau:M_B(X_1,n_1)\times M_B(X_2,n_2) \to M_B(X_1\times X_2, n_1n_2)$$
given by $\tau(V_1,V_2) = p_1^*V_1\otimes p_2^*V_2$. Also let
$$\sigma:M_B(X,n_1)\times M_B(X,n_2)\to M_B(X, n_1+n_2)$$
be the morphism $\sigma(V_1,V_2)\to V_1\oplus V_2$.
When the $X_i$ 
are smooth projective,  similar morphisms can be defined for $M_{Dol}$.

\begin{prop}\label{prop:kunneth}
 $\Sigma_{1,B}^k(X_1\times X_2\times \ldots X_N,n)$ is
the union of images of products
$$\sigma(\tau(\Sigma_{1,B}^{k_1}(X_1,n_1)\times\ldots
\Sigma_{1,B}^{k_N}(X_N,n_N))\times M_B(X_1\times\ldots
X_N,n-n_1n_2\ldots n_N))
$$
as  $k_i$ and $n_i$  range over partitions of $k$ and factorizations
of integers no greater than $n$ respectively.
\end{prop}

We will prove the proposition when $N=2$. Let
 $G_i =\pi_1(X_i)$, then of course $G=G_1\times G_2$ is the
 fundamental group of the produce. The classification of irreducible
 $G$-modules is well known, at least for finite groups.

\begin{lemma} Any irreducible finite dimensional complex representation
of $G$ is the form $V_1\otimes_\C V_2$ where $V_i$ are irreducible
$G_i$ representations.
\end{lemma}

\begin{proof} Let $V$ be a nonzero irreducible representation of $G$. Then
it contains a nonzero irreducible $G_1$-module $V_1$. The vector
space $G_2V_1 = \sum_{g\in G_2} gV_1$ is a $G$-module, so it
must coincide with $V$. Therefore $V$ is isomorphic to direct sum of
copies of $V_1$, so by standard arguments \cite[section 27]{cr},
 $R=End_{G_1}(V)$ is isomorphic to
the $\C$-endomorphism ring of an irreducible $R$-submodule $V_2$. $V_2$ is
necessarily a $G_2$-module, and the natural map $V_1\otimes V_2\to V$
is easily seen to be an isomorphism.  $V_2$ must be irreducible, since
otherwise, we could find a smaller submodule $V_2'$, and the $G$-module
$V_1\otimes V_2'$ would violate irreducibility of $V$.
\end{proof}

The above proposition is now a consequence of the usual K\"unneth
formula for cohomology. If $X$ is a product of subvarieties then,
there is a similar decomposition for $\Sigma_{Dol}^k$.

\section{Examples}

In this final  section, we will give some explicit examples.
To simplify notation, we will drop the decorations ``an, Dol, B'';
 it should  be clear from context where we are.
 In addition to the 
previous results, we will need the following well known fact:

\begin{lemma}
 The Euler characteristic of a rank $n$ local system $V$ on a variety $X$
is $n$ times the Euler characteristic $e(X)$ of $X$.
\end{lemma}
\begin{proof} 
There are a number of ways to 
prove this; perhaps the easiest is choose a finite triangulation of
$X$, then observe that each term $S^\dt(V)$ in the simplicial
cochain complex with coefficients in $V$ is (noncanonically) isomorphic to
$S^\dt(\C)^{\oplus n}$.
\end{proof}

\begin{ex} Let $C$ be a smooth projective curve of genus $g > 1$.
If $V$ is a nontrivial irreducible rank $n$ local system then 
$H^0(C,V) = 0$ and $H^2(C,V) \cong H^0(C,V^*)^*=0$.
 Therefore $dim\,H^1(C,V) = 2n(1-g)$. So that
$\Sigma_{2n(1-g),}^1(C,n) = M(C,n)$ because irreducible 
local systems are dense in $M(C,n)$. If $m > m_0=n(1-g)$, then
$\Sigma_{2m}^1(C,n)$ is locus of semisimple local systems with a trivial
summand of rank $m-m_0$. If $X= C_1\times C_2$ is a product of two
curves, then 
$$M(X,n) = 
\bigcup_{a+b=n}M(C_1,a)\times M(C_2,b) $$ and
$V\in \Sigma_1^1(X,n)$ if and only if it has a direct summand of
the form $p_i^*V_i$
 by  \ref{prop:kunneth}. Hitchin's map  respects this 
 decomposition. More precisely, the product
 $$h_{C_1}\times h_{C_2}:M(C_1,a)\times M(C_2,b) \to S_a(C_1)\times 
 S_b(C_2)$$
 factors through the restriction of $h_X$ to the above component, and
 the map 
 $$h_X( M(C_1,a)\times M(C_2,b))\to S_a(C_1)\times S_b(C_2)$$
 is  a bijection. To see this, let 
 $$(E,\theta)=(p_1^*E_1\otimes p_2^*E_2,
 p_1^*\theta_1\otimes I_{E_2} + I_{E_1}\otimes p_2^*\theta_2)$$
  be  a point
 on this component. 
 Then the projection of $h_X(E,\theta)$ to $S_a(C_1)\times S_b(C_2)$ is
 $(h_{C_1}(E_1,\theta_1), h_{C_2}(E_2,\theta_2))$.
 Moreover, the value of  $h_X(E,\theta)$ at any point 
 $x=(x_1,x_2)\in X$ is determined by the  eigenvalues of  
 $\theta_i(x_i)$ which is, in turn,  determined by 
 $h(E_i,\theta_i)(x_i)$.
\end{ex}

The generic vanishing theorem, shows that $\Sigma_1^k(X,n) $ is
proper whenever the $X$ has $1$-form without zeros. When $X$
is an abelian variety, we can say much more.

\begin{ex} Let $A$ be an abelian variety. Then for each $k$, 
$\Sigma_{1}^k(A,n)$ is locus of semisimple local systems with 
a trivial direct summand. To prove this, notice that 
$A$ is homeomorphic to a product of circles, and by \ref{prop:kunneth}, it
suffices to treat the case of a single circle. In this case its
clear, as $H^0(S^1,V)$ and $H^1(S^1,V)$ are respectively
just the invariants and coinvariants of $V$.
\end{ex}
 
\begin{ex}
Let $X\subset A$ be an ample 
divisors. Then $\Sigma_{1}^k(X,n)$ has the same description as in the 
previous example when $k < dim X$ by \ref{cor:lef}. 
The same holds for $k > dim X$ by duality \ref{cor:duality}. 
An easy Chern class argument shows that the Euler
characteristic of $X$ is nonzero. Consequently $\Sigma_{1}^{dim X}(X,n)$
is all of $M(X,n)$.
\end{ex}

The last example has an abelian fundamental group.
 Now we will show that any  group $\Gamma$ which 
can be realized as the fundamental group of a smooth projective
variety, is realizable as the fundamental group of smooth
projective variety  with $\Sigma_m^k(x,n)$ as large as possible, for
any $m,n$ and $k\ge 2$. Note that $\Sigma^1$ cannot be changed without 
altering the fundamental group.

\begin{ex} 
 Let $Y$ be smooth projective variety with fundamental group
$\Gamma$. We can arrange that $dim Y = k+1$ as follows. First after
replacing $Y$ with a product with a projective space, we can assume
that $dim Y \ge k+1$. Now replace $Y$ by a $k+1$ dimensional 
generic complete intersection.
Choose a very ample line bundle $O_Y(1)$ on $Y$, and let $X_d$ be the
general member of the linear system associated to $O_Y(d)$. Then
$\Sigma_m^i(Y) = \Sigma_m^i(X_d)= \Sigma_m^{2k-i}(X_d)$
 for $i < k$ by \ref{cor:lef} and
duality. Thus these sets are independent of $d$.
The Euler characteristic is just the Chern number $c_k(X_d)$,
and this can be computed to obtain a polynomial of degree $k$ in $d$. Thus
the Euler characteristic approaches $\pm\infty$ as $d\to \infty$.
Therefore after choosing $m$, we can find $d$
large enough so that $\Sigma_m^k(X_d,n) = M(X_d,n)$ 

\end{ex}


\begin{thebibliography}{MM}
\bibitem[A1]{a1}
D. Arapura. Higgs line bundles, Green-Lazarsfeld sets and maps of K\"ahler
manifolds to curves.  {\em Bull. A.M.S.} {\bf 26} (1992), 310-314.

\bibitem[A2]{a2}
D. Arapura. Geometry of cohomology support loci for local systems I.
{\em J. Alg. Geom.} to appear.

\bibitem[B]{b}
I. Biswas. A remark on the deformation theory of Green and Lazarsfeld.
{\em J. reine. ang. Math.} {\bf 449} (1994), 103-124

\bibitem[CR]{cr}
C. Curtis, I. Reiner, Representation theory of finite groups
and associative algebras, John Wiley \& Sons (1962)

\bibitem[D]{d}
S. Donaldson. Infinite determinants, stable bundles and curvature.
{\em Duke Math. J.} {\bf 54} (1987), 231-247

\bibitem[F]{fu}
A. Fujiki.  Hyperk\"ahler structure on the moduli space of flat
bundles. {\em in  Lect. Notes in Math 1468}, Springer-Verlag (1991)


\bibitem[GL1]{gl1}
M. Green and R. Lazarsfeld.  Deformation theory, generic vanishing theorems,
and some conjectures of Enriques, Catanese and Beauville. {\em Invent. Math.} 
{\bf 90} (1987), 389-407.

\bibitem[GL2]{gl2}
M. Green and R. Lazarsfeld.  Higher obstructions to
deforming cohomology groups of line bundles
{\em J. A.M.S.} {\bf 4} (1991), 87-103

\bibitem[EGA III]{ega3}
A. Grothendieck, J. Dieudonn\'e, El\'ements de g\'eom\'etrie
alg\'ebrique, {\em Publ. Math. I.H.E.S.} {\bf 11,17}

\bibitem[GM]{gm}
W. Goldman, J. Millson. The deformation theory of representations of 
the fundamental groupp of compact K\"ahler manifolds. {\em Publ. Math.
I.H.E.S.} {\bf 67} (1988), 43-96

\bibitem[GS]{gs}
 V. Guillemin, S. Sternberg, Symplectic techniques in physics,
Cambridge U. Press (1984)


\bibitem[H1]{hit1}
N. Hitchin.  The self-duality equations on a Riemann surface. 
{\em Proc.
London Math. Soc.}  {\bf 55} (1987), 59-126.

\bibitem[H2]{hit2}
N. Hitchin. Stable bundles and integrable systems. {\em Duke Math. J.} {\bf
54} (1987), 91-114.

\bibitem[H3]{hit3}
N Hitchin. Hyper K\"ahler manifolds, {\em Sem. Bourbaki}, (1992)

\bibitem[M]{mi}
J. Milnor. Morse theory. Princeton U. Press


\bibitem[N]{n}
 R. Narasimhan, On homology groups of Stein
spaces, {\em Invent. Math.} (1967)

\bibitem[P]{p}
T. Pantev. Lecture at Santa Cruz (1995)


\bibitem[S1]{s1}
C. Simpson.  Higgs bundles and local systems. {\em Publ. Math. I.H.E.S.} {\bf
75} (1992), 5-95.

\bibitem[S2]{s2}
C. Simpson.  Moduli of representations of the fundamental group of a smooth
projective variety, I: {\em Publ. Math. I.H.E.S.} {\bf 79} (1994), 47-129;
II: {\em Publ. Math. I.H.E.S.} {\bf 80} (1994), 5-79.

\bibitem[S3]{s3}
C. Simpson. Subspaces of moduli spaces of rank one local systems.  {\em Ann.
Sci. Ec. Norm. Sup.} {\bf 26} (1993), 361-401.

\bibitem[S4]{s4}
C. Simpson. The Hodge filtration in nonabelian cohomology.
{\em  alg-geom preprint} (1996)

\bibitem[Sp]{sp}
E. Spanier. Algebraic topology. McGraw-Hill (1966)

\bibitem[UY]{uy}
K. Uhlenbeck, S.T. Yau. On the existence of Hermitean-Yang-Mills 
connections in stable vector bundles, {\em Comm Pure and App. Math.}
{\bf 39} (1986),  257-293

\bibitem[V1]{v1}
M. Verbitsky, Deformation of trianalytic subvarieties of hyperk\"ahler
manifolds, {\em alg-geom preprint} (1996) 

\bibitem[V2]{v2}
M. Verbitsky, Desingularization of singular hyperk\"ahler manifolds 
I, II {\em alg-geom preprint} (1996)



\end{thebibliography}
\end{document}